%% file: ArXivBruchhaeuser.tex
\documentclass[a4paper,10pt]{article}

\usepackage[utf8]{inputenc}
\usepackage{amsfonts}
\usepackage{amssymb,amsmath,amsthm}

\usepackage{authblk}

\usepackage{tikz}
\usetikzlibrary{plotmarks}
\usetikzlibrary{shapes,arrows}
\usetikzlibrary{positioning}
\usetikzlibrary{trees,mindmap}
\usetikzlibrary{fadings,shapes.arrows,shadows}
\usetikzlibrary{decorations.pathreplacing}
\usetikzlibrary{calc}
\usetikzlibrary{patterns}
\usepackage{pgfplots}
\pgfplotsset{compat=1.10}
\usepgfplotslibrary{fillbetween}

\usepackage{array}
\usepackage{hhline}

\usepackage{booktabs}

\usepackage{subfig}

\usepackage{graphicx}

\newtheorem{defi}{Definition}[section]
\newtheorem{theorem}[defi]{Theorem}

\newtheorem{remark}[defi]{Remark}

\numberwithin{equation}{section}
\numberwithin{table}{section}
\numberwithin{figure}{section}

\begin{document}

\title{\Large \textbf{Dual weighted residual based error control for nonstationary 
convection-dominated equations: potential or ballast?}}
\author[M.\ P.\ Bruchh\"auser, M.\ Bause]
{\large  M.\ P.\ Bruchh\"auser\thanks{bruchhaeuser@hsu-hh.de}\,, 
K.\ Schwegler,
M.\ Bause\thanks{bause@hsu-hh.de}\\
{\small Helmut Schmidt University, Faculty of Mechanical Engineering,\\ 
Holstenhofweg 85, 22043 Hamburg, Germany}
}
\date{}
\maketitle

\begin{abstract}
Even though substantial progress has been made in the numerical approximation of 
convection-dominated problems, its major challenges remain in the scope of current 
research \cite{BruchhaeuserJKN18}. In particular, parameter robust a posteriori error 
estimates for quantities of physical interest and adaptive mesh refinement strategies 
with proved 
convergence are still missing. Here, we study numerically the potential of the Dual 
Weighted Residual (DWR) approach applied to stabilized finite element methods to further 
enhance the quality of approximations. The impact of a strict application of the DWR 
methodology is particularly focused rather than the reduction of computational costs for 
solving the dual problem by interpolation or localization.
\end{abstract}


\section{Introduction}
\label{BruchhaeuserSec:Intro}
In their recent review paper \cite{BruchhaeuserJKN18} the authors nicely survey the 
current state of research in the numerical approximation of convection-dominated 
equations and incompressible flow. These problems with prominent applications in many 
branches of technology have strongly attracted researchers' interest not just 
since the pioneering works of the 1980's (cf., e.g., \cite{BruchhaeuserBH81,BruchhaeuserHB79}). 
The introduction of various families of linear and nonlinear residual-based 
stabilization techniques and of algebraic stabilization techniques, usually 
referred to as algebraic flux correction schemes, are regarded as milestones in 
the development of discretization techniques that are able to reduce spurious 
nonphysical oscillations close to sharp layers of convection-dominated problems;
cf. \cite{BruchhaeuserJS08} for a comparative study of those techniques.
For a general review of those methods as well as a list of references we refer to, 
e.g., \cite{BruchhaeuserJKN18,BruchhaeuserRST08}. In \cite{BruchhaeuserJKN18}, 
the authors further identify numerous problems that are still unresolved in this field of 
research. In particular, the non-availabilty of parameter-robust a posteriori error 
estimates for quantities of physical interest and in general situations is 
stressed. Moreover, the authors point out that adaptive mesh refinement strategies 
that are based on such a posteriori error estimates and guarantee convergence in 
appropriate norms are desirable and indespensable for further improvement.

One possible technique for those adaptive strategies is the DWR 
method (\cite{BruchhaeuserBR98,BruchhaeuserBR01,BruchhaeuserBR03}), where the
error is estimated in an arbitrary user-chosen goal quantity of physical interest.
The DWR approach relies on a space-time variational formulation of the discrete 
problem and uses duality techniques to find rigorous a posteriori error estimates,
obtained through the approximation of an additional dual problem.
Early studies for adaptive mesh refinement applied to various stationary stabilized equations 
date back to the end of the last century; cf. \cite[Section 3.3 and 8]{BruchhaeuserBR01}
for a brief overview and further literature. For the nonstationary Navier-Stokes 
equations the DWR approach was applied together with local projection stabilization 
(LPS) in \cite{BruchhaeuserBR12}. In \cite{BruchhaeuserAJ15} an adaptive algorithm
in time is presented for convection-dominated problems, where the time step control
uses a post-processed solution.

In this work, we study numerically the potential of combining the DWR approach with 
SUPG (\cite{BruchhaeuserHB79,BruchhaeuserBH81}) stabilized finite element methods for the 
efficient and reliable approximation of nonstationary convection-dominated problems. 
Here, a \emph{first dualize and then stabilize} principle is applied; cf. 
Sect.~\ref{BruchhaeuserSec:3}.
For the approximation of the dual problem higher-order schemes are used, due to
the recently received results for stationary convection-dominated problems in
\cite{BruchhaeuserBSB18}, where a comparative study to an approximation by 
higher-order interpolation was done.
The presented numerical results illustrate the performance properties and robustness
of the proposed algorithm with respect to a vanishing perturbation or diffusion 
parameter. Thereby, the \emph{potential} of a DWR approach is demonstrated
and it is shown that the DWR based adaptivity is not \emph{ballast} for the approximation
of convection-dominated problems. 

This work is organized as follows. In Sect.~\ref{BruchhaeuserSec:2} we present
our model problem as well as the stabilized space-time discretization.
A first dualize and then stabilize DWR approach as well as a localized error 
representation is presented in Sect.~\ref{BruchhaeuserSec:3}. Finally, in  
Sect.~\ref{BruchhaeuserSec:numerical} the results of our numerical experiments are
presented.

\section{Model problem and stabilized space-time discretizations}
\label{BruchhaeuserSec:2}
In this work we consider the following convection-diffusion-reaction
equation 
\begin{equation}
\begin{array}{rclcl}
\partial_{t} u-\nabla \cdot (\varepsilon \nabla 
 u) + \vec{b}\cdot \nabla u + \alpha\,u
 & = & f & \mathrm{in } & \Omega\times I\,,  \\
 u & = & g_D & \mathrm{on } & 
 \partial\Omega\times I\,,\\
 u(0) & = & u_{0} & \mathrm{in } & 
 \Omega\,,
\end{array}
\label{Bruchhaeuser:cdr}
 \end{equation}
where $\Omega \subset \mathbf{R}^{d}$, with $d=2$ or $d=3$ is a polygonal or polyhedral 
bounded domain with Lipschitz boundary $\partial\Omega$ and $I=(0, T), T>0$, 
is a bounded domain in time. 
To ensure the well-posedness of problem (\ref{Bruchhaeuser:cdr}) we assume that
$0 < \varepsilon \leq 1$ is a constant diffusion coefficient, 
$\vec{b} \in \vec{H}^1(\Omega)\cap\vec{L}^{\infty}(\Omega)$ is the flow field or convection 
tensor, $\alpha \in L^{\infty}(\Omega)$ is the reaction coefficient,
$u_0\in H_0^1(\Omega)$ is a given initial condition, 
$f\in L^{2}(I;L^2(\Omega))$ is a given source of the unknown scalar 
quantity $u$ and $g\in L^{2}(I;H^{\frac{1}{2}}(\partial\Omega))$
is a given function specifying the Dirichlet boundary condition.
Furthermore, we assume that the conditions 
$\nabla \cdot \vec{b}(\vec{x})=0$ and  $\alpha(\vec{x})\geq 0$ 
are fulfilled for all $\vec{x}\in\Omega$.
Henceforth, for the sake of simplicity, we deal with homogeneous Dirichlet boundary 
conditions only. In our numerical examples in Sect. 
\ref{BruchhaeuserSec:numerical}, we also consider more general boundary conditions;
cf. also Remark \ref{BruchhaeuserRem:InhomDirichletBC}. 

It is well known that problem (\ref{Bruchhaeuser:cdr}) along with the above conditions
admits a unique weak solution
$u \in V:=\big\{v\in L^{2}(I; H_0^1(\Omega))\big|
\;\partial_{t}v\in L^{2}(I; H^{-1}(\Omega))\big\}$ that satisfies the following 
variational formulation; cf., e.g. \cite{BruchhaeuserRST08}.

\textit{Find $u \in V$, satisfying $u(0)=u_0$, such that}
\begin{equation}
 A(u)(\varphi) = F(\varphi) \quad \forall \varphi \in V\,,
 \label{Bruchhaeuser:weak}
\end{equation}
\textit{\indent where the bilinear form} $A:V\times V \to \mathbf{R}$ 
\textit{and the linear form} $F:V\to \mathbf{R}$ \textit{are}
\begin{eqnarray}
 \label{Bruchhaeuser:bilinearform_A}
 A(u)(\varphi) & := & \sum_{n=1}^{N}\int_{I_n}\big\{(\partial_t u,\varphi)
 +a(u)(\varphi)
 \big\} \mathrm{d} t 
 + \sum_{n=2}^{N}([u]_{n-1},\varphi_{n-1}^+ )\,, \\
\,F(\varphi) & := & \int_I(f,\varphi)\;\mathrm{d}t\,.
 \label{Bruchhaeuser:linearform_F}
\end{eqnarray}
Here, $a(u)(\varphi):=(\varepsilon\nabla u, \nabla \varphi) +(\vec{b}\cdot \nabla u, \varphi) + (\alpha u,\varphi)$
and $(\cdot,\cdot)$ denotes the standard inner product of $L^2(\Omega)$.
For the discretization in time  we divide the time interval $I$ into not necessarily 
equidistant, left-open subintervals $I_n:=(t_{n-1},t_n]\,,$ with $n=1,\dots,N\,,$ where 
$0=t_0<t_1<\dots<t_N=T$ with step size $\tau_n=t_n-t_{n-1}$ and 
$\tau=\max_n\,\tau_{n}\,.$ Next, we introduce the time-discrete 
function spaces.
\begin{eqnarray}
 V_{\tau}^{cG(r)} & := & \Big\{ v\in C(\bar{I}; H_0^1(\Omega))\big|\; v|_{I_{n}}\in 
 \mathbf{P}_{r}(\bar{I}_{n}; H_0^1(\Omega)) 
 \Big\}\,.
 \label{Bruchhaeuser:time_spaces_cG}\\
 V_{\tau}^{dG(r)} & := & \Big\{ v\in L^{2}(I; H_0^1(\Omega))\big|\; v|_{I_{n}}\in 
 \mathbf{P}_{r}(I_{n}; H_0^1(\Omega)),
  v_\tau(0)\in L^2(\Omega)\Big\}\,,
 \label{Bruchhaeuser:time_spaces_dG}
 \end{eqnarray}
where $\mathbf{P}_{r}(\bar{I}_{n}; H_0^1(\Omega))$ denotes the space of all 
polynomials in time up to degree $r\geq0$ on $I_n$ with values in $H_0^1(\Omega)\,.$
For some function $v_\tau\in V_{\tau}^{dG(r)}$ we define the limits $v_\tau^\pm$ from
above and below of $v_\tau$ as well as their jump at $t_n$ by
$$v_{\tau,n}^\pm := \lim_{t\to t_n\pm0} v_\tau(t) \,, \qquad [v_\tau]_n := v_{\tau,n}^+ 
-v_{\tau,n}^- \,.$$
Using the discontinuous Galerkin method for the time
discretization of the so called primal problem (\ref{Bruchhaeuser:weak}) leads to the 
following time-discrete variational approximation.

\textit{Find $u_\tau \in V_{\tau}^{dG(r)}$ such that}
\begin{equation}
 A(u_\tau)(\varphi_\tau) + (u_{\tau,0}^+,\varphi_{\tau,0}^+) = 
 F(\varphi_\tau) + (u_{0},\varphi_{\tau,0}^+) \quad \forall \varphi_\tau \in V_{\tau}^{dG(r)}\,,
 \label{Bruchhaeuser:time_discrete_scheme}
\end{equation}
\textit{\indent with} $A(\cdot)(\cdot)$ \textit{and} $F(\cdot)$ \textit{being defined by} (\ref{Bruchhaeuser:bilinearform_A})
\textit{and} (\ref{Bruchhaeuser:linearform_F}), \textit{respectively}.\\
We note that the initial condition is incorporated into the variational problem.
Next, we describe the Galerkin finite element approximation in space of the semidiscrete
problem (\ref{Bruchhaeuser:time_discrete_scheme}).
We use Lagrange type finite element spaces of continuous functions that are piecewise 
polynomials. For the discretization in space, we consider a decomposition 
$\mathcal{T}_{h}$
of the domain $\Omega$ into disjoint elements $K$, such that 
$\bar{\Omega}=\cup_{K\in\mathcal{T}_{h}}\bar{K}$. Here, we choose the elements 
$K\in\mathcal{T}_{h}$ to be quadrilaterals for $d=2$ and hexahedrals for $d=3$. 
We denote by $h_{K}$ the diameter of the element $K$. The global space 
discretization parameter $h$ is given by $h:=\max_{K\in\mathcal{T}_{h}}h_{K}$. 
Our mesh adaptation process yields locally refined cells, which is 
enabled by using hanging nodes. We point out that 
the global conformity of the finite element approach is preserved since the unknowns 
at such hanging nodes are eliminated by interpolation between the neighboring 
'regular' nodes; cf.~\cite{BruchhaeuserBR03}.
On $\mathcal{T}_{h}$ we define the discrete finite element space by
$ V_{h}^{p,n}:=
  \big\{v\in V\cap C(\bar{\Omega})\,\big|\,v|_{K}
  \in Q_h^p(K)
  \,,\forall K\in\mathcal{T}_{h}, \big\}\,,$
$n=1,\dots,N $, where $Q_h^p(K)$ is the space of polynomials that are of degree less
than or equal to $p$ with respect to each variable $x_1,\dots,x_d$.
By replacing $H_0^1(\Omega)$ in the definition of the semidiscrete function 
spaces $V_{\tau}^{cG(r)}$ and $V_{\tau}^{dG(r)}$ in (\ref{Bruchhaeuser:time_spaces_cG}) and 
(\ref{Bruchhaeuser:time_spaces_dG}), respectively, by $V_h^{p,n}$, we obtain the fully discrete 
function spaces
\begin{eqnarray}
V_{\tau h}^{cG(r),p} & := & \left\{v_{\tau h}\in V_\tau^{cG(r)} \big| v_{\vert 
I_n}\in \mathbf{P}_r(\bar{I}_n;V_h^{p,n})
\right\}
\label{Def:Xkhrp_cG}
\\
V_{\tau h}^{dG(r),p} & := & \left\{v_{\tau h}\in V_\tau^{dG(r)} \big| v_{\vert 
I_n}\in \mathbf{P}_r(I_n;V_h^{p,n})
\,, v_{\tau h}(0) \in V_h^{p}\right\} 
\label{Def:Xkhrp_dG}
\end{eqnarray}
with $V_{\tau h}^{cG(r),p}\subset V_\tau^{cG(r)}$ and 
$V_{\tau h}^{dG(r),p}\subset V_\tau^{dG(r)}$. We note that the spatial finite 
element space $V_h^{p,n}$ is allowed to be different on 
all intervals $I_n$ which is natural in the context of a discontinuous Galerkin 
approximation of the time variable and allows dynamic mesh changes in time.
The fully discrete discontinuous in time scheme then reads as follows.

\textit{Find $u_{\tau h} \in V_{\tau h}^{dG(r),p}$ such that}
\begin{equation}
 A(u_{\tau h})(\varphi_{\tau h}) + (u_{\tau h,0}^+,\varphi_{\tau h,0}^+) = 
 F(\varphi_{\tau h}) + (u_{0},\varphi_{\tau h,0}^+) \quad \forall \varphi_{\tau h} \in V_{\tau h}^{dG(r),p}\,,
 \label{Bruchhaeuser:space_time_discrete_scheme}
\end{equation}
\textit{\indent with} $A(\cdot)(\cdot)$ \textit{and} $F(\cdot)$ \textit{being defined in} (\ref{Bruchhaeuser:bilinearform_A})
\textit{and} (\ref{Bruchhaeuser:linearform_F}), \textit{respectively}. \\
In this work we focus on convection-dominated problems with small difussion 
parameter $\varepsilon$. Then the finite element approximation needs to be 
stabilized in order to reduce spurious and non-physical oscillations of the 
discrete solution arising close to sharp fronts or layers. 
Here, we apply the streamline upwind Petrov-Galerkin method (for short
SUPG); cf. \cite{BruchhaeuserHB79,BruchhaeuserBH81,BruchhaeuserRST08,BruchhaeuserJS08,BruchhaeuserBS12}. 
The stabilized variant of the fully discrete scheme 
(\ref{Bruchhaeuser:space_time_discrete_scheme}) then reads as follows.

\textit{Find $u_{\tau h} \in V_{\tau h}^{dG(r),p}$ such that}
\begin{equation}
 A_{S}(u_{\tau h})(\varphi_{\tau h}) + (u_{\tau h,0}^+,\varphi_{\tau h,0}^+) = 
 F(\varphi_{\tau h}) + (u_{0},\varphi_{\tau h,0}^+) \quad \forall \varphi_{\tau h} \in V_{\tau h}^{dG(r),p}\,,
 \label{Bruchhaeuser:stabilized_space_time_discrete_scheme}
\end{equation}
\textit{\indent with} $A_{S}(u)(\varphi):=A(u)(\varphi)+S(u)(\varphi)$ 
\textit{and stabilization terms}
\begin{displaymath}
\begin{array}{l@{}l}
S(u_{\tau h})(\varphi_{\tau h}) &{} :=  \displaystyle\sum_{n=1}^N\int_{I_n} 
\sum\limits_{K\in \mathcal{T}_h}\delta_K\big( \mathit{R}(u_{\tau h}), 
\vec{b}\cdot\nabla\varphi_{\tau h}\big)_K \,\mathrm{d} t \\
+ \displaystyle\sum\limits_{n = 2}^{N}\sum\limits_{K\in\mathcal{T}_h} &{}\hspace{-0.2cm} \delta_K
\big(\left[u_{\tau h}\right]_{n-1}, 
\vec{b}\cdot\nabla\varphi_{\tau h,n-1}^+\big)_{K} + \displaystyle\sum\limits_{K\in\mathcal{T}_h}
\delta_K \big( u_{\tau h,0}^{+} - u_0,
\vec{b}\cdot\nabla\varphi_{\tau h,0}^{+} \big)_{K}\,,\vspace{0.2cm} \\
\hspace{0.8cm}\mathit{R}(u_{\tau h}) &{} :=  \partial_{t} u_{\tau h} + \vec{b}\cdot\nabla u_{\tau h} 
- \nabla\left(\varepsilon\nabla u_{\tau h}\right)  + \alpha u_{\tau h} - f\,.
\end{array}
\end{displaymath}
\begin{remark}
The proper choice of the stabilization parameter $\delta_K$ is an
important issue in the application of the SUPG approach; 
cf., e.g.,~\cite{BruchhaeuserJN11} and the discussion 
therein. For the situation of steady-state convection and reaction,
an optimal error estimate for $\delta_K=\mathrm{O}(h)$ was derived in 
\cite{BruchhaeuserJN11}.
\end{remark}

\section{A DWR approach and a posteriori error estimation}
\label{BruchhaeuserSec:3}
The DWR method aims at the control of an error in an arbitrary user-chosen target 
functional $J$ of physical interest. To get an error representation with respect to this 
target functional, an additional dual problem has to be solved.
Before we focus on the error representation, we introduce the
dual problem of (\ref{Bruchhaeuser:weak}) whose derivation is based on the 
Euler-Lagrangian method of constrained optimization. 
For a detailed derivation we refer to \cite[Chapter 6,9]{BruchhaeuserBR03}.
We note that here a so called first dualize and then stabilize principle
is used, where the stabilization is applied to the discrete dual problem 
after its derivation via the Euler-Lagrangian method of constrained optimization;
cf. \cite[Remark 2]{BruchhaeuserBSB18}.

For some given functional 
$J:V \to \mathbf{R}$ we consider solving
\begin{displaymath}
J(u) = \min\{J(v)\,, \; v\in V\,, \; \mathrm{where}\; A(v)(\varphi)=F(\varphi)\; \forall 
\varphi \in V\}\,.
\end{displaymath}
We assume that the functional $J$ is Fr\'{e}chet differentiable. i.e. $J'(y)\in V'$
for $y\in V$.
For the derivation of the error representation we define the corresponding 
Lagrangian functional $\mathcal{L}:V \times V \to \mathbf{R}$ by
\begin{equation}
 \mathcal{L}(u,z):=J(u)+F(z)-A(u)(z)-(u(0)-u_{0},z(0))\,,
 \label{BruchhaeuserDef:L}
\end{equation}
where we refer to $z\in V$ as the dual variable (or Lagrangian multiplier); 
cf.~\cite{BruchhaeuserBR03}. We determine a stationary point $\{u,z\}\in V\times V$ of 
$\mathcal{L}(\cdot,\cdot)$ by the condition 
$\mathcal{L}'(u,z)(\psi,\varphi)\\=0$,
or equivalently by the system of equations 
\begin{eqnarray*}
A'(u)(\psi,z) & = & J'(u)(\psi) \quad \forall \psi \in V\,, 
 \\[1ex]
A(u)(\varphi) & = & F(\varphi) \quad \quad \;\; \forall 
\varphi \in V\,.
\end{eqnarray*}
The second of these equations, the $z$-component of the stationary condition, is just 
the given primal problem (\ref{Bruchhaeuser:weak}),
whereas the $u$-component of the stationary  condition, is called the 
dual or adjoint equation with $A'(u)(\psi,z) = \int_I\big\{(\partial_t \psi,z)
 +a(\psi, z) \big\} \mathrm{d}t$ and $J'(u)(\psi)=\int_I\big\{(j(u),\psi)\big\} \mathrm{d}t$
 for some function $j(\cdot)\in L^2(I;L^2(\Omega))$.
Applying integration by parts in time to the first term of $A'$ and taking the condition
$\nabla \cdot \vec{b}(\vec{x})=0$ into account (cf. Sect. \ref{BruchhaeuserSec:2}) 
yields the representation
$A^\ast(z)(\psi):= 
A'(u)(\psi,z) 
= \int_I\big\{-(\partial_t z,\psi)
 +(\varepsilon\nabla z, \nabla \psi)-(\vec{b}\cdot\nabla z,\psi)+(\alpha z,\psi) \big\} \mathrm{d}t\,.$
Finally, we find by using the proposed stabilized Galerkin discretization scheme 
(\ref{Bruchhaeuser:stabilized_space_time_discrete_scheme}) the following stabilized discrete dual
problem. 

\textit{Find $z_{\tau h} \in V_{\tau h}^{dG(r),p}$ such that}
\begin{equation}
 A_{S}^{\ast}(z_{\tau h})(\psi_{\tau h}) + (z_{\tau h,T}^-,\psi_{\tau h,T}^-) = 
 J'(u_{\tau h})(\psi_{\tau h}) \quad \forall \psi_{\tau h} \in V_{\tau h}^{dG(r),p}\,.
 \label{Bruchhaeuser:dual_space_time_discrete_scheme}
\end{equation}
In (\ref{Bruchhaeuser:dual_space_time_discrete_scheme}), we put 
$
A_{S}^\ast(z_{\tau h})(\psi_{\tau h}) := A^\ast (z_{\tau h})(\psi_{\tau h}) 
+ S^\ast(z_{\tau h})(\psi_{\tau h})
$
with
\begin{displaymath}
\label{Bruchhaeuser:As}
\begin{array}{r@{}l}
A^\ast(z_{\tau h})(\psi_{\tau h}) := &{} \displaystyle \sum_{n=1}^{N}  \int_{I_n} 
\big\{-(\partial_t z_{\tau h}, \psi_{\tau h})- (\vec{b}\cdot \nabla z_{\tau h},\psi_{\tau h})
+ (\varepsilon \nabla z_{\tau h}, \nabla \psi_{\tau h})\big\} \mathrm{d} t  \\
&{} + (\alpha z_{\tau h}, \psi_{\tau h}) - \displaystyle\sum\limits_{n=2}^{N} 
 (\left[z_{\tau h}\right]_{n-1},\psi_{\tau h,n-1}^{-})\,, \\
S^\ast(z_{\tau h})(\psi_{\tau h}) := &{} \displaystyle\sum_{n=1}^{N}\int_{I_n} 
\sum\limits_{K\in\mathcal{T}_h}\delta_K^\ast  
(\mathit{R}^\ast(z_{\tau h}), -\vec{b}\cdot\nabla\psi_{\tau h})_K \mathrm{d} t\\[1ex]
- \displaystyle\sum\limits_{n = 2}^{N}\sum\limits_{K\in\mathcal{T}_h}
&{}\delta_K^\ast\big(\left[z_{\tau h}\right]_{n-1}, -\vec{b}\cdot\nabla\psi_{\tau h,n-1}^{-}\big)_{K} 
 + \displaystyle\sum\limits_{K\in\mathcal{T}_h}
\delta_K^\ast 
(z_{\tau h,N}^{-} ,- \vec{b}\cdot\nabla\psi_{\tau h,N}^- )_{K}\,,\\[2ex]
\mathit{R}^\ast(z_{\tau h}) := &{} -\partial_{t} z_{\tau h} -\vec{b}\cdot\nabla z_{\tau h} 
- \nabla\left(\varepsilon\nabla z_{\tau h}\right)  + \alpha z_{\tau h} - j(u_{\tau h})\,.
\end{array}
\end{displaymath}
To derive a representation of the error $J(e)=J(u) - J(u_{\tau h})$ we 
need some abstract results. 
In order to keep this work self-contained we pare down to the key arguments
of the DWR approach applied to the stabilized model problem. We follow 
the lines of \cite[Chapter 6 and 9]{BruchhaeuserBR03} and \cite{BruchhaeuserB00}, where
all of the proofs can be found.
To start with, we need to extend the definition of the Lagrangian functional to 
arguments of $(V + V_{\tau h}^{dG(r),p})\times V$. 
In the following we let $ \mathcal{L}:(V+ V_{\tau h}^{dG(r),p})\times V$ 
be defined by 
\begin{equation}
\label{Bruchhaeuser:Lext}
\mathcal{L} (u,z) := 
J(u)+F(z)- A(u)(z)
- \displaystyle\sum\limits_{n=2}^{N}  
\big(\left[u\right]_{n-1},z_{n-1}^+\big) - \big(u(0)-u_0,z(0)\big)\,.
\end{equation}
Then it follows that 
\begin{equation}
\label{Bruchhaeuser:SC}
\mathcal{L}_{u} (u,z)(\psi) + \mathcal{L}_{z}(u,z)(\varphi) = 0 \quad \forall
\{\psi,\varphi\}\in V\times V\,.
\end{equation}
The discrete solution $\{u_{\tau h},z_{\tau h}\} \in V_{\tau h}^{dG(r),p}\times 
V_{\tau h}^{dG(r),p}$ then satisfies
\begin{equation}
\label{Bruchhaeuser:DSSC}
\mathcal{L}_{u} (u_{\tau h}, z_{\tau h})(\psi_{\tau h}) 
+  \mathcal{L}_{z}(u_{\tau h},z_{\tau h})(\varphi_{\tau h})
 =  S(u_{kh})(\varphi_{kh}) + S^\ast(z_{kh})(\psi_{kh})
\end{equation}
for all
$\{\psi_{\tau h},\varphi_{\tau h}\} \in V_{\tau h}^{dG(r),p}\times V_{\tau h}^{dG(r),p}$.
For the defect of the discrete solution in the stationary condition 
(\ref{Bruchhaeuser:DSSC}) we use the notation 
\begin{displaymath}
\tilde{S}(x_{\tau h})(y_{\tau h}) := S(u_{\tau h})(\varphi_{\tau h}) 
+ S^{\ast}(z_{\tau h})(\psi_{\tau h})\,,
\end{displaymath}
with $x_{\tau h}:= \{u_{\tau h},z_{\tau h}\} \in V_{\tau h}^{dG(r),p}\times 
V_{\tau h}^{dG(r),p}$ and $y_{\tau h}:=\{\psi_{\tau h},\varphi_{\tau h}\}\in 
V_{\tau h}^{dG(r),p}\times V_{\tau h}^{dG(r),p}$. 
To derive a representation of the error $J(u) - J(u_{\tau h})$ 
we need the following abstract theorem that develops the error in terms of the 
Lagrangian functional.
\begin{theorem}
 Let $X$ be a function space and $\mathcal{L}: X \rightarrow \mathbf{R}\;$ be a
 three times differentiable functional on $X$. Suppose that $x_c\in X_c$ with 
 some (''continuous'') function space 
$X_c \subset X$ is a stationary point of $\mathcal L$. Suppose that $x_d \in X_d$ with 
some (''discrete'') function space $X_d \subset X\,,$ with not necessarily $X_d \subset 
X_c\,,$ is a Galerkin approximation to $x_c$ being defined by the equation 
\begin{displaymath}
\mathcal{L}'(x_d)(y_d) = \tilde{S}(x_d)(y_d) \quad  \forall y_d \in X_d\,.
\end{displaymath}
In addition, suppose that the auxiliary condition $\mathcal{L}'(x_c)(x_d) = 0$ 
is satisfied. Then there holds the error representation
\begin{displaymath}
\mathcal{L}(x_c) - \mathcal{L}(x_d) = \frac{1}{2} \mathcal{L}'(x_d) (x_c- y_d) + 
\frac{1}{2} 
\tilde{S}(x_d)(y_d -x_d) + \mathcal{R} \quad \forall y_d \in X_d\,,
\end{displaymath}
where the remainder $\mathcal{R}$ is defined by
$\mathcal{R} = \frac{1}{2}\int_0^1 \mathcal{L}'''(x_d +s e)(e,e,e)
\cdot s \cdot (s-1)\,\mathrm{d} s\,,$
with the notation $e:=x_c-x_d$.
\label{BruchhaeuserThm:L}
\end{theorem}
For the subsequent theorem we introduce the primal and dual residuals by
\begin{eqnarray}
\rho(u_{\tau h})(\varphi) & := &  \;F(\varphi)-A(u_{\tau h})(\varphi) 
-(u_{\tau h,0}^{+}-u_{0},\varphi(0))
 \quad\;\;\;\;  \forall \varphi \in V\,, 
\label{BruchhaeuserDef:primalresidual}\\
\rho^\ast(z_{\tau h})(\psi) & := &  \;J'(u_{\tau h})(\psi)-A^\ast(z_{\tau h})(\psi)   
-(z_{\tau h,N}^{-},\psi(T))
 \quad  \forall \psi \in V\,.
\label{BruchhaeuserDef:dualresidual}
\end{eqnarray}
\begin{theorem}
\label{BruchhaeuserThm:J}
Suppose that $\{u,z\}\in V \times V$ is a stationary point of 
the Lagrangian functional $\mathcal{L}$ defined in (\ref{Bruchhaeuser:Lext}) 
such that (\ref{Bruchhaeuser:SC}) is satisfied. 
Let $\{u_{\tau h},z_{\tau h}\}\in V_{\tau h}^{dG(r),p} \times V_{\tau h}^{dG(r),p}$ 
denote its Galerkin approximation being defined by 
(\ref{Bruchhaeuser:stabilized_space_time_discrete_scheme}) and 
(\ref{Bruchhaeuser:dual_space_time_discrete_scheme}), respectively,  
such that (\ref{Bruchhaeuser:DSSC}) is satisfied. 
Then there holds the error representation that 
\begin{equation}
J(u) - J(u_{\tau h}) = \frac{1}{2}\rho(u_{\tau h}) (z - \varphi_{\tau h}) + 
\frac{1}{2}\rho^\ast(z_{\tau h}) (u-\psi_{h}) + \mathcal{R}_{\tilde{S}} + 
\mathcal{R}_{J}
\label{BruchhaeuserThm:J1}
\end{equation}
for arbitrary functions $\{\varphi_{\tau h},\psi_{\tau h}\}\in V_{\tau h}^{dG(r),p} 
\times V_{\tau h}^{dG(r),p}$, where 
the remainder terms are \
$\mathcal{R}_{\tilde{S}} := \frac{1}{2} S(u_{\tau h})(\varphi_{\tau h}+z_{\tau h}) + \frac{1}{2} 
S^\ast(z_{\tau h})(\psi_{\tau h}-u_{\tau h})$
and 
$\mathcal{R}_{J} := \frac{1}{2}\int_0^1 J'''(u_{\tau h}+s\cdot~e)(e,e,e) 
\cdot s\cdot (s-1)\,\mathrm{d} s\,,$
with $e=u-u_{\tau h}$.
\end{theorem}
In the error respresentation (\ref{BruchhaeuserThm:J1}) the continuous solution 
$u$ is required for the evaluation of the dual residual. The following theorem  
shows the equivalence of the primal and dual residual up to a quadratic remainder. 
This observation will be used below to find our final error respresentation in 
terms of the goal quantity $J$ and a suitable linearization for its computational 
evaluation or approximation, respectively.  
\begin{theorem}
\label{BruchhaeuserThm:ResDev}
Under the assumptions of Thm.~\ref{BruchhaeuserThm:J}, and with the definitions 
(\ref{BruchhaeuserDef:primalresidual}) and (\ref{BruchhaeuserDef:dualresidual}) 
of the primal and dual residual, respectively, there holds that 
\begin{displaymath}
\rho^\ast(z_{\tau h}) (u-\psi_{\tau h}) = \rho(u_{\tau h}) (z - \varphi_{\tau h}) + 
S(u_{\tau h})(\varphi_{\tau h}-z_{\tau h}) 
+ S^\ast (z_{\tau h})(u_{\tau h}-\psi_{\tau h}) + \Delta \rho_{J}\,,
\end{displaymath}
for all $\{\psi_{\tau h},\varphi_{\tau h}\}\in V_{\tau h}^{dG(r),p} 
\times V_{\tau h}^{dG(r),p}$, where the remainder term is 
given by \\
$\Delta \rho_{J}:=  - \int_0^1 J''(u_{\tau h} + s \cdot e)(e,e) \,\mathrm{d} s$
with $e:=u - u_{\tau h}$.
\end{theorem}
In a final step we combine the results of the previous two theorems to get a 
localized approximation of the error that is then used for the design of the 
adaptive algorithm.
We note that the final result (\ref{BruchhaeuserEq:localER}) is a
slight modification of Thm. 5.2 for the nonstationary Navier-Stokes equations
stabilized by LPS in \cite{BruchhaeuserBR12}. 
The difference comes through using a first dualize and then stabilize approach 
as well as SUPG stabilization here.
\newpage
\begin{theorem}[Localized error representation]
\label{BruchhaeuserThm:LocPresent}
Let the assumptions of Thm.~\ref{BruchhaeuserThm:J} be satisfied.
Neglecting the higher-order error terms, 
then there holds as a linear approximation the cell-wise error representation that
\begin{equation}
\label{BruchhaeuserEq:localER}
\begin{array}{r@{}l}
J(u) - J(u_{\tau h})  \doteq &{} \displaystyle \sum_{n=1}^{N}\int_{I_n}\sum\limits_{K\in\mathcal{T}_h} \Big\{ 
\big(\mathit{R}(u_{\tau h}), z-\varphi_{\tau h}\big)_{K} 
- \delta_K\big(\mathit{R}(u_{\tau h}), \mathbf{b}\cdot\nabla\varphi_{\tau h}\big)_{K} \\[1ex]
&{} - \big(\mathit{E}(u_{\tau h}), z-\varphi_{\tau h}\big)_{\partial K}\Big\}\mathrm{d}t \\[1ex]
&{} -\displaystyle \sum_{K\in\mathcal{T}_{h}}
\left(u_{\tau h,0}^{+}-u_{0},z(t_{0})-\varphi_{\tau h,0}^{+}\right)_{K}\\[1ex]
&{} - \displaystyle \sum_{n=2}^{N}\sum_{K\in\mathcal{T}_{h}}
\left([u_{\tau h}]_{n-1},
z(t_{n-1})-\varphi_{\tau h,n-1}^{+}\right)_K\\[1ex]
&{} + \displaystyle \sum_{K\in\mathcal{T}_{h}}\delta_{K}\left(u_{\tau h,0}^{+}-u_{0},
\vec{b}\cdot\nabla\varphi_{\tau h,0}^{+}\right)_K \\[1ex]
&{} + \displaystyle \sum_{n=2}^{N}\sum_{K\in\mathcal{T}_{h}}\delta_{K}\left(
[u_{\tau h}]_{n-1},\vec{b}\cdot\nabla\varphi_{\tau h,n-1}^{+}
\right)_K
\,.
\end{array}
\end{equation}
The cell- and edge-wise residuals are defined by 
\begin{eqnarray}
\label{Bruchhaeusereq:42} \mathit{R}(u_{\tau h})_{|K} &:=& f -\partial_t u_{\tau h}
+ \nabla\cdot(\varepsilon\nabla 
u_{\tau h})
- \mathbf{b}\cdot\nabla u_{\tau h} -\alpha u_{\tau h} \,,\\[0.5ex]
\label{Bruchhaeusereq:43} \mathit{E}(u_{\tau h})_{|\Gamma} &:=& \left\{ \begin{array}{cl} 
\frac{1}{2}\mathbf{n}\cdot[\varepsilon\nabla u_{\tau h}] & 
\mbox{ if } \Gamma\subset\partial K\backslash\partial\Omega\,, \\[0.5ex] 
0 & 
\mbox{ if } \Gamma\subset\partial\Omega\,,\\ \end{array}\right.
\end{eqnarray}
where $[\nabla u_{\tau h}]:= \nabla u_{\tau h}{}_{|\Gamma\cap K}
-\nabla u_{vh}{}_{|\Gamma\cap K'}$ defines the jump of $\nabla u_{\tau h}$ 
over the inner edges $\Gamma$ with normal unit vector 
$\mathbf{n}$ pointing from $K$ to $K'$. 
\end{theorem}
\begin{remark}
 \label{BruchhaeuserRem:InhomDirichletBC} (Nonhomogeneous Dirichlet boundary 
conditions)
In the case of nonhomogeneous Dirichlet boundary conditions the following 
additional term has to be added to the error representation (\ref{BruchhaeuserEq:localER})
\begin{eqnarray*}
\displaystyle \sum_{n=1}^{N}\int_{I_n}
-\big((g_D-\tilde{g}_{D,\tau h}),\varepsilon\nabla z\cdot\mathbf{n}\big)_{\partial\Omega}
\mathrm{d}t\,,
\end{eqnarray*}
where the discrete function $\tilde{g}_{D,\tau h}$ is an appropriate finite element 
approximation of the extension $\tilde{g}_D$ in the sense that the trace 
of $\tilde{g}_D$ equals $g_D$ on $\partial\Omega$; cf. 
\cite{BruchhaeuserBSB18,BruchhaeuserBR03}.
\end{remark}

\section{Numerical studies}
\label{BruchhaeuserSec:numerical}
In this section we illustrate and investigate the performance properties of the 
proposed approach of combining the DWR method with SUPG stabilized finite element
approximations of nonstationary convection-dominated problems.
Therefore some general indications are needed. The error representation
(\ref{BruchhaeuserEq:localER}), written as
\begin{equation}
 J(u) - J(u_{\tau h})  \doteq \eta:= \displaystyle \sum_{n=1}^{N}
 \sum\limits_{K\in\mathcal{T}_h} \eta_{K}^{n}\,, 
 \label{BruchhaeuserEq:DefEta}
\end{equation}
depends on the discrete primal solution $u_{\tau h}$ as well as on the exact 
dual solution $z$. For solving the primal problem (\ref{Bruchhaeuser:cdr})
we use the discontinuous in time scheme 
(\ref{Bruchhaeuser:stabilized_space_time_discrete_scheme})
to get a discrete solution $u_{\tau h} \in V_{\tau h}^{dG(r),p}$. For the 
application  of (\ref{BruchhaeuserEq:DefEta}) in computations, the 
unknown dual solution~$z$ has to be approximated, which results in an approximate 
error indicator $\tilde{\eta}$. 
This approximation cannot be done in the same finite element space as used for 
the primal problem, since this would result in an useless vanishing error 
representation $\tilde{\eta}=0$, due to Galerkin orthogonality. In contrast to the
approximation by higher-order interpolation which is widespread used in the literature,
cf. \cite{BruchhaeuserBR03,BruchhaeuserB00,BruchhaeuserBR12},
we use an approximation by higher-order finite elements here. This is done due to
the results in \cite{BruchhaeuserBSB18}, where a comparative study between these two
approaches is presented for steady convection-dominated problems. In this study 
the superiority of using higher-order finite elements was shown for an increasing 
convection dominance.
Thus, we use for the discretization of the dual problem a finite element space
that consists of polynomials in space and time that are at least of one polynomial 
degree higher than its primal counterpart, more precisely we compute a discrete
dual solution $z_{\tau h}\in V_{\tau h}^{cG(r+1),p+1}$.
For the now following example we briefly present our algorithm. For further details 
we refer to \cite{BruchhaeuserKBB17}.

\noindent\rule{\textwidth}{1pt}
  \begin{center}
   \textbf{Adaptive solution algorithm (Refining in space and time)}
  \end{center}
\vspace{-0.3cm}
\noindent\rule{\textwidth}{0.5pt}
\textbf{Initialization:} Set $i=0$ and generate the initial space-time slab
$Q=\Omega\times I$ with $Q:=\cup_{n=1}^{N}Q_n=\cup_{n=1}^{N}(\Omega\times I_n)$.
\begin{enumerate}
   \item Compute the \textbf{primal} and \textbf{dual} solution 
   $u_{\tau h} \in V_{\tau h}^{dG(r),p}$ and $z_{\tau h} \in V_{\tau h}^{cG(r+1),p+1}$
   \item Evaluate the \textbf{a posteriori space-time error indicator}
  $ \tilde{\eta}  :=  \sum_{n=1}^{N}\sum_{K\in\mathcal{T}_h} \tilde{\eta}_K^n\,.\, $ \\
  Mark the time intervals $I_{\tilde{n}}$ where $I_{\tilde{n}}$ belongs to the set of the
  time intervals $I_n$ according to $\theta_\tau$ percent of the 
  worst indicators $\tilde{\eta}^n:=\sum_{K\in\mathcal{T}_h} |\tilde{\eta}_K^n|$.
  Mark those cells $\tilde{K}$ (of the respective spatial mesh of a space-time slab
  $Q_n$)
  for refinement  that make up a certain fraction $\theta_h$ of the total error.
  \item Check the \textbf{stopping condition}:
  If  $\tilde{\eta} < \textrm{tol}$ is satisfied, then the  
  adaptive solution algorithm is terminated.
  \item Else, \textbf{adapt} the space-time slab $Q$,
   increase $i$ to $i+1$ and return to Step~1.
  \end{enumerate}
  \vspace{-0.3cm}
\noindent\rule{\textwidth}{0.5pt}

For the implementation of the adaptive algorithm we use our \texttt{DTM++} 
frontend software \cite{BruchhaeuserKBB17} that is based on the open 
source finite element library \texttt{deal.II}; cf.~\cite{BruchhaeuserAABBBDGHHKKMPTW}.
For measuring the accuracy of the error estimator, we will 
study 
the effectivity index 
\begin{equation}
\label{BruchhaeuserEq:Ieff}
\mathcal{I}_{\mathrm{eff}} = \left|\frac{\tilde{\eta}}{J(u)-J(u_{h})}\right|
\end{equation}
as the ratio of the estimated error $\tilde{\eta}$ of 
(\ref{BruchhaeuserEq:DefEta}) over the exact error. 
Desirably, the index $\mathcal{I}_{\mathrm{eff}}$ should be close to one.

\subsection{Example  (Rotating hill with changing orientation).}
\label{BruchhaeuserExample1}

\begin{figure}[!ht]
\centering
\includegraphics[trim=4.5cm 19.8cm 4.5cm 3.6cm,width=\textwidth]{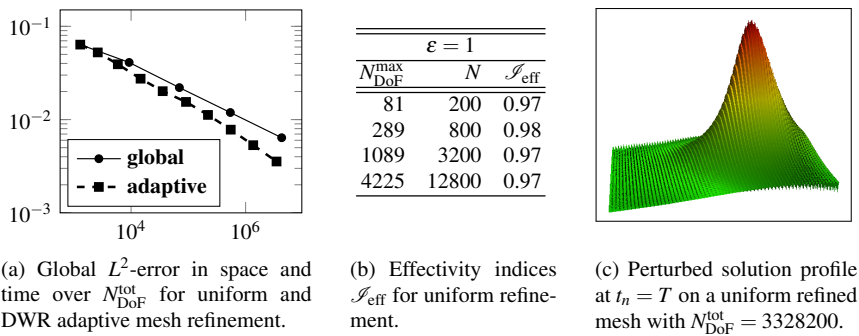}
\caption{Comparison of global space-time $L^2$-errors, table for uniform refinement
as well as perturbed solution profile.}
\label{BruchhaeuserFig:KB1andIeffandSolProf}
\end{figure}

In this example we analyze the performance properties and the robustness of our 
algorithm with respect to the small perturbation parameter $\varepsilon$. 
We study problem (\ref{Bruchhaeuser:cdr}) with the prescribed solution
\begin{equation}
\label{BruchhaeuserEq:movghump}
u(t,x,y):= \frac{\arctan(5\pi(2t-1))}{1+a_0\big(x-\frac{1}{2}-\frac{1}{4}\cos(2\pi t)\big)^2+
 a_0\big(y-\frac{1}{2}-\frac{1}{4}\sin(2\pi t)\big)^2}\,,
\end{equation}
where $\Omega\times I := (0,1)^2\times (0,1]$ and $a_0 = 50$. 
We choose the flow field $\mathbf{b} = (2,3)^\top$ and the reaction coefficient $\alpha = 
1.0$. 
The solution (\ref{BruchhaeuserEq:movghump}) is characterized by a counterclockwise
rotating hill and designed in a way that the orientation of the hill changes its sign
from negative to positive at the midpoint of the time interval $I$ at $t=0.5$.
For the solution (\ref{BruchhaeuserEq:movghump}) the right-hand side function $f$ is 
calculated from the partial differential equation. Boundary conditions 
are given by the exact solution. Our target quantity is chosen to control the 
global $L^2$-error in space and time, given by
\begin{equation}
\label{BruchhaeuserEq:HumpTarget}
J(u)= \frac{1}{\|e\|_{(0,T)\times\Omega}}\displaystyle\int_I(u,e)\mathrm{d}t\,,
\quad \mathrm{with} \;\; \|\cdot\|_{(0,T)\times\Omega}=\left(\int_I(\cdot,\cdot)\;\mathrm{d}t\right)^{\frac{1}{2}}\,.
\end{equation}
In our first test we investigate problem (\ref{Bruchhaeuser:cdr}) for $\varepsilon=1$
and without any stabilization to verify our algorithm for a simple test case.
In Fig.~\ref{BruchhaeuserFig:KB1andIeffandSolProf}b we monitor the effectivity indices for 
uniform refinement.
Here and in the following, $N_{\mathrm{DoF}}^{\mathrm{tot}}$ denotes the total 
number of degrees of freedom in space and time for one DWR loop while
$N_{\mathrm{DoF}}^{\mathrm{max}}$ denotes the maximum number 
of degrees of freedom of a spatial mesh  used 
within one DWR loop. Furthermore, $N$ denotes the total number of space-time slabs $Q_n$
used for one DWR loop; 
cf. the adaptive solution algorithm at the beginning of this chapter.
As we are using continuous finite elements of lowest order in space and
discontinuous finite elements of lowest order in time (for short cG(1)-dG(0)), we
refine globally once in space and twice in time after each step. 
For this uniform refinement the algorithm provides convincing and reliable
values for the related effectivity indices; 
cf. the table in Fig. \ref{BruchhaeuserFig:KB1andIeffandSolProf}.

In the sequel, all investigations are performed using SUPG stabilization
in the sense of (\ref{Bruchhaeuser:stabilized_space_time_discrete_scheme}) for
varying diffusion coefficients. Initially, in Fig.~\ref{BruchhaeuserFig:KB1andIeffandSolProf}a 
we compare the convergence behavior of the proposed DWR approach with a uniform 
mesh refinement strategy for $\varepsilon=10^{-6}$. The DWR based 
adaptive mesh adaptation is clearly superior to the uniform refinement in terms 
of accuray over the total number of degrees of freedom in space and time.
\begin{table} 
\centering
\begin{minipage}{\linewidth}
\centering
\resizebox{0.98\linewidth}{!}{%
\begin{tabular}{r@{\quad}r@{\quad}r@{\quad}c@{\quad}c@{\quad}c@{\quad}|r@{\quad}r@{\quad}r@{\quad}c@{\quad}c@{\quad}c@{\quad}}
\toprule
\toprule
\multicolumn{6}{c@{\quad}|}{$\varepsilon=10^{-3}$} & \multicolumn{6}{c}{$\varepsilon=10^{-6}$}
\\
\midrule
$N_{\mathrm{DoF}}^{\mathrm{tot}}$ & $N_{\mathrm{DoF}}^{\mathrm{max}}$ & $N$ & $J(e)$ & $|\tilde{\eta}|$ & $\mathcal I_{\mathrm{eff}}$ &
$N_{\mathrm{DoF}}^{\mathrm{tot}}$ & $N_{\mathrm{DoF}}^{\mathrm{max}}$ & $N$ & $J(e)$ & $|\tilde{\eta}|$ & $\mathcal I_{\mathrm{eff}}$
\\
\midrule
17650    & 628   & 36  & 3.244e-02 & 2.301e-02 &  0.71 & 10974   & 617   & 22  & 4.761e-02 & 2.335e-02 & 0.491\\ 
37728    & 961   & 54  & 2.381e-02 & 9.258e-03 &  0.39 & 22868   & 951   & 33  & 3.317e-02 & 1.824e-02 & 0.550\\ 
86919    & 1464  & 81  & 1.839e-02 & 1.667e-02 &  0.91 & 52537   & 1454  & 49  & 2.507e-02 & 1.883e-02 & 0.751\\ 
196074   & 2310  & 121 & 1.286e-02 & 4.921e-03 &  0.38 & 121978  & 2357  & 73  & 2.039e-02 & 3.193e-02 & 1.566\\ 
409573   & 3744  & 181 & 8.943e-03 & 7.127e-03 &  0.80 & 264748  & 3718  & 109 & 1.378e-02 & 5.095e-03 & 0.370\\ 
942465   & 6196  & 271 & 6.017e-03 & 4.767e-03 &  0.79 & 631452  & 5786  & 163 & 9.734e-03 & 8.891e-03 & 0.913\\ 
2135099  & 10326 & 406 & 4.753e-03 & 3.706e-03 &  0.78 & 1506529 & 9258  & 244 & 7.409e-03 & 4.234e-03 & 0.572\\ 
4678474  & 15715 & 609 & 3.101e-03 & 3.128e-03 &  1.01 & 3610055 & 15167 & 366 & 4.772e-03 & 3.673e-03 & 0.770\\ 
10179407 & 24309 & 913 & 2.231e-03 & 2.240e-03 &  1.00 & 8109271 & 25747 & 549 & 3.314e-03 & 3.374e-03  & 1.018\\ 
\bottomrule
\end{tabular}}
\end{minipage}
\caption{Effectivity indices for the goal quantity (\ref{BruchhaeuserEq:HumpTarget}) 
for different values of $\varepsilon$.}
\label{BruchhaeuserTab:KB1Ieff}
\end{table}
In Table \ref{BruchhaeuserTab:KB1Ieff} we present selected effectivity indices 
of the proposed DWR approach applied to the stabilized approximation scheme 
(\ref{Bruchhaeuser:stabilized_space_time_discrete_scheme}) for different
diffusion coefficients.    
For $\varepsilon=10^{-3}$ the effectivity indices nicely converge to 
one for an increasing number of the total degrees of freedom in space and time.
For the more challenging case of $\varepsilon=10^{-6}$ the values of 
$\mathcal I_{\mathrm{eff}}$ are still close to one, which confirms the robustness
of our adaptive algorithm with respect to the small perturbation parameter 
$\varepsilon$.

\begin{figure}[!ht]
\centering
\includegraphics[trim=4.2cm 18.7cm 4.5cm 3.3cm,width=\textwidth]{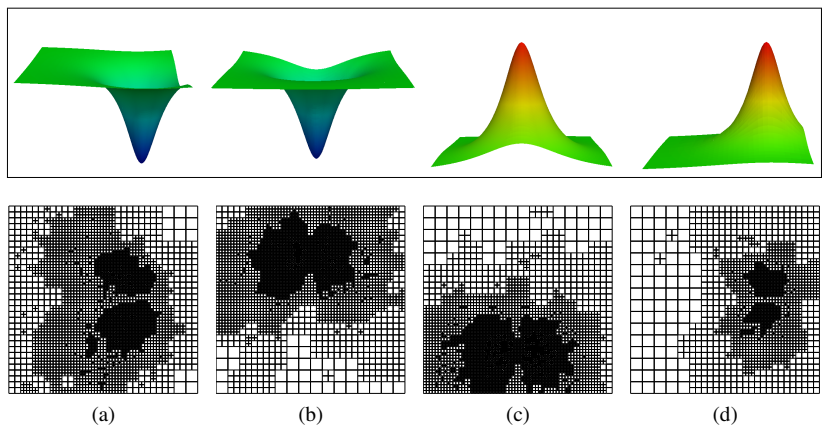}
\caption{Stabilized solution profiles and related adaptive spatial meshes after 
9 DWR loops at time points 
$t_n=0$ (a), $t_n=0.25$ (b),
$t_n=0.75$ (c) and $t_n=1$ (d) 
for $\varepsilon=10^{-6}$.}
\label{BruchhaeuserFig:SolutionProfiles}
\end{figure}

To finish this section, we illustrate in Fig.~\ref{BruchhaeuserFig:SolutionProfiles}
some stabilized solution profiles as well as adaptive spatial meshes after the 
last DWR loop at selected time points. Comparing the stabilized solution profile
obtained on an adaptive refined mesh (cf.~Fig.~\ref{BruchhaeuserFig:SolutionProfiles}d) with its 
unstabilized counterpart obtained on an uniform refined mesh 
(cf.~Fig.~\ref{BruchhaeuserFig:KB1andIeffandSolProf}c), both computed at the final 
time point $T$, we point out that the occuring spurious oscillations are strongly
reduced. Furthermore, no smearing effects are observed, comparing the solution profiles 
of the hill at the respective time points. Considering the underlying
spatial meshes, we note that the total number of the spatial cells is more or less
equal during the whole time period, but the arrangement differs depending on
the related position of the hill. Thus, for the chosen target functional
(\ref{BruchhaeuserEq:HumpTarget}) the spatial mesh runs as expected with the 
rotation of the hill in a synchronous way. In addition, we note that the 
mesh refinement is slightly weaker at the final time point. 
This is due to the error propagation of the underlying problem which is captured 
by the dual weights in the error estimate. This effect is in good agreement to
the results obtained for the heat equation in \cite[p. 122]{BruchhaeuserBR03}.

\section{Conclusions and Outlook}
\label{BruchhaeuserSec:summary}
In this work we presented an adaptive solution algorithm for SUPG stabilized
finite element approximations of time-dependent convection-diffusion-reaction
equations. The underlying approach is based on the Dual Weighted Residual method 
for goal-oriented error control. A \emph{first dualize and then stabilize} 
philosophy was applied for combining the space-time adaptation process in the 
course of the DWR approach with the stabilization of the finite element techniques. 
We used a higher-order finite element approximation in space and time
in order to compute the dual solution. 
In numerical experiments we could prove that spurious oscillations that 
typically arise in numerical approximations of convection-dominated problems 
could be reduced significantly. Robust effectivity indices close to one were 
obtained for different values of the diffusion coefficient. 
It was shown that the DWR based adaptivity is no \emph{ballast} on the way to solve 
convection-dominated problems. Conversely, it offers \emph{potential} for further 
improvements in handling those problems. 

We note that recent results in post-processing variational time 
discretization schemes (cf., e.g., \cite{BruchhaeuserBKRS18,BruchhaeuserAJ15}) allow 
the computation of improved solutions admitting an additional order of convergence 
for the discretization in time by negligible computational costs, and thus offer 
further potential for reducing the costs of computing the dual solution.

\subsection*{Acknowledgement}
We acknowledge Uwe K\"ocher for his help in the design and implementation
of the underlying software DTM++/dwr-diffusion; cf.
\cite{BruchhaeuserKBB17}.

\input{ArXivBruchhaeuser_referenc}

\end{document}

%% file: ArXivBruchhaeuser_referenc.tex
%
%
%